\documentclass[11pt]{amsart}
\usepackage[T1]{fontenc}
\usepackage{textcomp}
\usepackage[margin=1in]{geometry}
\usepackage{mathtools}
\usepackage{amssymb}
\usepackage{enumitem}
\usepackage{xcolor}

\usepackage{listings}
\definecolor{codedarkgreen}{RGB}{51, 133, 4}
\definecolor{codemaroon}{RGB}{133, 5, 63}
\definecolor{codeteal}{RGB}{0, 145, 109}
\definecolor{codepurple}{RGB}{123, 35, 125}
\newcommand{\GRevLex}{\textcolor{codeteal}{\texttt{GRevLex}}}
\newcommand{\Lex}{\textcolor{codeteal}{\texttt{Lex}}}
\newcommand{\Limit}{\textcolor{codeteal}{\texttt{Limit}}}
\newcommand{\true}{\textcolor{codeteal}{\texttt{true}}}
\newcommand{\false}{\textcolor{codeteal}{\texttt{false}}}
\newcommand{\cache}{\textcolor{codeteal}{\texttt{cache}}}
\newcommand{\gens}{\textcolor{blue}{\texttt{gens}}}
\newcommand{\generators}{\textcolor{blue}{\texttt{generators}}}
\newcommand{\Ideal}{\textcolor{codedarkgreen}{\texttt{Ideal}}}
\newcommand{\HashTable}{\textcolor{codedarkgreen}{\texttt{HashTable}}}
\newcommand{\GroebnerBasis}{\textcolor{codedarkgreen}{\texttt{GroebnerBasis}}}
\lstdefinelanguage{Macaulay2}{
basicstyle=\normalsize\ttfamily,
  alsoletter=",
  classoffset=1,
  keywords={transpose,det,subsets,genericMatrix,needsPackage,presentation,generators,gens,selectInSubring,i1,i2,i3,i4,i5,i6,i7,i8,i9,i10,i11,flatten,ideal},
  keywordstyle={\color{blue}},
classoffset=2,
breaklines=true,
morekeywords={"SubalgebraBases","Engine"},
keywordstyle={\color{codemaroon}},
classoffset=3,
morekeywords={QQ,CacheTable,Matrix,Eliminate},
keywordstyle={\color{codedarkgreen}},
classoffset=4,
morekeywords={restart,false,true,Weights,Limit,Lex,MonomialOrder},
keywordstyle={\color{codeteal}},
classoffset=5,
morekeywords={list,for,in,from,to,of},
keywordstyle={\color{codepurple}},
xleftmargin=1em,
xrightmargin=1em,
columns=fullflexible,
keepspaces=true,
stepnumber=1,
numbers=none,
captionpos=b,
showspaces=false,
frame=none
}

\newcommand{\quoteoption}[1]{\textcolor{codemaroon}{\texttt{#1}}}
\newcommand{\QQ}{\mathbb{Q}}
\newcommand{\RR}{\mathbb{R}}

\usepackage{xcolor}
\usepackage{hyperref}
\DeclareMathOperator{\initial}{in}
\DeclareMathOperator{\SE}{SE}
\DeclareMathOperator{\SO}{SO}
\DeclareMathOperator{\se}{\mathfrak{s}\mathfrak{e}}

\theoremstyle{definition}
\newtheorem{definition}{Definition}[section]
\newtheorem{example}[definition]{Example}
\newtheorem{remark}[definition]{Remark}

\newcommand{\bC}{{\mathbb{C}}}
\newcommand{\bQ}{{\mathbb{Q}}}
\newcommand{\bP}{\mathbb P}
\usepackage{amssymb} 
\def\acts{\curvearrowright}
\newcommand{\subduction}{\texttt{subduction}}
\newcommand{\subring}{{\texttt{Subring}}}
\newcommand{\Macaulay}{{\textit{Macaulay2}}} 
\newcommand{\SubalgberaBases}{{\texttt{SubalgebraBases}}}
\newcommand{\SAGBIBasis}{{\texttt{SAGBIBasis}}}
\title{SubalgebraBases in Macaulay2}
\author[Burr]{Michael Burr}
\address{220 Parkway Drive, Clemson University, Clemson, SC 29634}
\email{burr2@clemson.edu}
\thanks{Burr was partially supported by NSF grant DMS-1913119 and SIMONS collaboration grant \#964285.}
\author[Clarke]{Oliver Clarke}
\address{The University of Edinburgh, James Clerk Maxwell Building, Edinburgh EH9 3FD} 
\thanks{The work was started when Oliver Clarke was at the University of Bristol supported by the EPSRC Doctoral Training Partnership
(DTP) award EP/N509619/1, and continued while being an overseas researcher under Postdoctoral Fellowship of Japan Society for the Promotion of Science (JSPS).}
\email{oliver.clarke@ed.ac.com}
\author[Duff]{Timothy Duff}
\address{Department of Mathematics, University of Washington, Seattle, WA 98195}
\thanks{Duff acknowledges support from an NSF Mathematical Sciences Postdoctoral Research Fellowship (DMS-2103310)}
\email{timduff@uw.edu}
\author[Leaman]{Jackson Leaman}
\address{220 Parkway Drive, Clemson University, Clemson, SC 29634}
\email{jleaman@g.clemson.edu}
\author[Nichols]{Nathan Nichols}
\address{Department of Mathematics and Statistical Sciences, Marquette University, Milwaukee, WI 53233}
\email{nathannichols454@gmail.com}

\author[Walker]{Elise Walker}
\address{Center for Computing Research, Sandia National Laboratories, Albuquerque, NM, USA}
\email{eawalke@sandia.gov}
\date{\today}

\begin{document}
\maketitle

\begin{abstract}
We describe a recently revived version of the software package \SubalgberaBases, which is distributed in the \Macaulay\ computer algebra system.
The package allows the user to compute and manipulate subagebra bases ---  which are also known as SAGBI bases or canonical bases and form a special class of Khovanskii bases --- for polynomial rings and their quotients.
We provide an overview of the design and functionality of \SubalgberaBases\ and demonstrate how the package works on several motivating examples.
\end{abstract}

\section{Introduction}

In the field of computational algebraic geometry, the notion of a Gr\"{o}bner basis plays an instrumental role in solving many basic problems involving polynomial ideals.
For computations involving polynomial {\em algebras}, there is an analogous, but more subtle, notion of a \emph{subalgebra basis}.  
This notion was introduced independently by Kapur and Madlener~\cite{KapurMadlener:1989} and Robbiano and Sweedler~\cite{RobbianoSweedler:1990}, where the names \emph{canonical basis} and \emph{SAGBI basis} were also proposed.
In more recent years, the name \emph{Khovanskii basis} has also been adopted~\cite{KM} to describe a more general notion for valued algebras.

The \SubalgberaBases\ package in \Macaulay\ was initially developed around 1997 by Mike Stillman and Harrison Tsai~\cite{M2}. 
In subsequent years this package was inaccessible to most users due to internal changes in \Macaulay, until version \texttt{1.18} was released in mid--2021.
Our work before and after this re-release has focused on restoring the package's original functionality, implementing additional algorithms, and designing new data structures that facilitate working with larger examples and admit new features in the package.

\medskip
\noindent\textbf{Outline.} In Section~\ref{sec:background}, we provide some background on subalgebra bases, with examples illustrating the basic usage of our package. 
In Section~\ref{sec:data-structures}, we introduce the main data structures: the \texttt{Subring} and \texttt{SAGBIBasis} types, which allow for resuming partial subalgebra basis computations,
and describe the various methods that they enable. 
Section~\ref{sec:options} summarizes several options that may be useful for computing subalgebra bases, and Section~\ref{sec:other-functionality}
covers some additional functionality.
Finally, in Section~\ref{sec: examples}, we provide some more sophisticated examples of our package in action. 
The examples appearing in this paper may be found in the file \texttt{accompanyingCode.m2}.

\section{Background and basic computations}\label{sec:background}

For a field $k$, let $R\vcentcolon=k[x_1,\dots,x_n]$ be a polynomial ring, which we equip with a monomial order $<$.
Unless otherwise specified, we follow the \Macaulay\ convention for variable ordering, in which $x_n < \dots  < x_1.$
Consider a subalgebra $A \vcentcolon= k[\{f_i\}|_{i\in I}] \subset R,$ where $I$ is some index set.
The {\em initial algebra} $\initial_<(A)$ is generated by the initial terms of all elements in $A$:
\[
\initial_< (A)=k[\initial_<(f):f\in A].
\]
\begin{definition}\label{def:SAGBI}
A set $\{g_j\}_{j\in J}\subseteq A$
is a {\em subalgebra basis} for $A$, with respect to the order $<$, if \[
\initial_<(A)=k[\{\initial_<(g_j)\}_{j\in J}].\]
\end{definition}

Just as Gr\"{o}bner bases enable many computations with polynomial ideals, the knowledge of a finite subalgebra basis allows us to answer several basic questions about a given algebra.
A basic application is \emph{algebra membership}: if $g_1, \ldots , g_s \in A$ form a subalgebra basis for $A,$ then any polynomial $f\in R$ has an associated \emph{normal form} $r\in R$: for some polynomial $q \in k[y_1, \ldots , y_s],$
\begin{equation}\label{eq:normal-form}
f = q (g_1, \ldots , g_s) + r,
\end{equation}
with the property that either $\initial_< (r) \notin \initial_< (A)$ or  $r=0.$
Additionally, we require that none of the monomials of $r$ are contained in $\initial_< (A)$, which implies it is unique.
The polynomials $q$ and normal form $r$ appearing in Equation~\eqref{eq:normal-form} can be computed using an analogue of multivariate polynomial division known as the \emph{subduction algorithm}. 
A description of the subduction algorithm may be found, for example, in~\cite[Algorithm 1.5]{RobbianoSweedler:1990} or~\cite[Algorithm 11.1]{Sturmfels:1996}.

Unlike Gr\"{o}bner bases, the theory of subalgebras bases is not, strictly speaking, algorithmic. 
For instance, as illustrated in Example~\ref{ex:degreeLimitedInfiniteSubalgebraBasis},
a finitely generated polynomial algebra need not have a finite subalgebra basis with respect to a given term order.
Nevertheless, there are many interesting examples for which finite subalgebra bases exist and can be computed.
We consider a classical and well-known example from invariant theory.


\begin{example}\label{ex:basicSubalgebra}
Let $R=\QQ [x_1,x_2,x_3]$ and consider the first three power-sums in $R$,
\begin{align*}
f_1 &= 
\underline{x_1}+x_2+x_3,\\
f_2 &= \underline{x_1^2}+x_2^2+x_3^2,\\
f_3 &= \underline{x_1^3}+x_2^3+x_3^3.
\end{align*}
The underlined monomials above are the initial terms with respect to the \GRevLex\ order on $R.$
The algebra $A = \QQ[f_1, f_2, f_3] $ is the ring of invariants for the standard action of the symmetric group $S_3 \acts R,$ given by $\sigma \cdot (x_1, x_2, x_3) = (x_{\sigma (1)}, x_{\sigma (2)}, x_{\sigma (3)}).$
A subalgebra basis for $A$ consists of the three elementary symmetric polynomials,
\begin{align*}
g_1 &= \underline{x_1}+x_2+x_3,\\
g_2 &= \underline{x_1x_2}+x_1x_3+x_2x_3,\\
g_3 &= \underline{x_1x_2x_3},
\end{align*}
hence $\initial_< (A) = \QQ [x_1, x_1 x_2, x_1 x_2 x_3].$ We verify these assertions in \Macaulay\ with the code below. 
\begin{lstlisting}[language=Macaulay2]
needsPackage "SubalgebraBases";
R = QQ[x_1..x_3];
A = subring {x_1+x_2+x_3, x_1^2+x_2^2+x_3^2, x_1^3+x_2^3+x_3^3};
SB = sagbi A
isSAGBI SB
\end{lstlisting}
The two lines without semicolons produce the following output, informing us that \texttt{SB} is an instance of the type \texttt{SAGBIbasis}, and that our computation terminated successfully.
\begin{lstlisting}[language=Macaulay2]
i4 : SB = sagbi A
o4 = SAGBIBasis Computation Object with 3 generators, Limit = 20.
o4 : SAGBIBasis
i5 : isSAGBI SB
o5 = true
\end{lstlisting}
To verify $f = x_1^4 + x_2^4 + x_3^4 \in A,$ we may compute $q$ and $r$ appearing in Equation~\eqref{eq:normal-form} as follows:
\begin{lstlisting}[language=Macaulay2]
i5 : A = subring(gens SB, GeneratorSymbol => g);
i6 : f = x_1^4 + x_2^4 + x_3^4;
i7 : q = f // A

      4     2       2
o7 = g  - 4g g  + 2g  + 4g g
      1     1 2     2     1 3
     
o7 : QQ[g ..g ]
         1   3
i8 : r = f % A
o8 = 0
o8 : R
\end{lstlisting}
\end{example}

For a given set of algebra generators $\{ f_i \}_{i\in I}$, the \emph{binomial syzygies} on the set of initial terms $\{ in_< (f_i) \}_{i\in I}$, also known as \textit{t\^ete-a-t\^etes}~\cite{RobbianoSweedler:1990}, provide an analogue of S-polynomials from Gr\"{o}bner basis theory.
The binomial syzygies $y^\alpha - y^\beta $ generate the kernel of the monomial map
\begin{align*}
\varphi_{in_< (f)} : k[y_i \mid i \in I] &\mapsto k[\initial_<(f_i)\mid i \in I],\\
y_i &\mapsto \initial_<(f_i).
\end{align*}
If we consider a presentation of the algebra $A$ given by
\begin{align*}
\varphi_{f} : k[y_i \mid i \in I] &\mapsto A,\\
y_i &\mapsto f_i,
\end{align*}
then the initial term of the ``S-polynomial'' $S_{\alpha , \beta} = \varphi_f (y^\alpha - y^\beta ) \in A$ is strictly less than $\initial_< \left( \varphi_f (y^\alpha) \right)=\initial_< \left( \varphi_f (y^\beta) \right)$.
Applying the subduction algorithm, we obtain a normal form $r_{\alpha, \beta}\in A$ for $S_{\alpha , \beta }$ analogous to $r$ in Equation~\eqref{eq:normal-form}.
If every $r_{\alpha , \beta }$ is zero, then the set $\{ f_i \}_{i \in I}$ forms a subalgebra basis.
Otherwise, we enlarge the set of generators to include all nonzero $r_{\alpha , \beta }$ and restart the computation.

\begin{example}\label{ex:degreeLimitedInfiniteSubalgebraBasis}
Consider the algebra $A = \QQ [x_1+x_2, x_1x_2, x_1x_2^2] \subseteq \QQ [x_1, x_2]$, and let $<$ be the \GRevLex\ order. The initial algebra $\initial_< (A) = \QQ [x_1, x_1 x_2, x_1x_2^2, x_1x_2^3, \dots]$ is not finitely generated. However, we can still compute a partial subalgebra basis in \Macaulay\ as follows:
\begin{lstlisting}[language=Macaulay2]
i3 : R = QQ[x_1,x_2];
i4 : SB = subalgebraBasis({x_1+x_2, x_1*x_2, x_1*x_2^2}, Limit => 7)
o4 = | x_1+x_2 x_1x_2 x_1x_2^2 x_1x_2^3 x_1x_2^4 x_1x_2^5 x_1x_2^6 |
             1       15
o4 : Matrix R  <--- R
\end{lstlisting}
The option \Limit\ (whose default value is 20) allows the procedure outlined above to terminate.
The final S-polynomial computed is $S_{\alpha , \beta} = (x_1+x_2)(x_1x_2^{5}) - (x_1x_2^2)(x_1x_2^{4}) = x_1 x_2^{6}$, of degree $7.$
Using \texttt{isSAGBI}, we can check that the computation is incomplete:
\begin{lstlisting}[language=Macaulay2]
i3 : isSAGBI SB
o3 = false
\end{lstlisting}
\end{example}

Various extensions of Definition~\ref{def:SAGBI} may be found in the literature.
For instance, in~\cite{Miller}, a 
Noetherian integral domain with identity replaces the coefficient field $k$. 
Building upon this, the authors of \cite{SS} extend the theory of subalgebra bases to quotients of these polynomial rings.  
Currently, our package can be used to compute a (partial) subalgebra basis for a finitely-generated subalgebra $A \subseteq R / I$, where $I$ is a polynomial ideal. 
In this setting, the initial algebra $\initial_< (A)$ is defined to be a subalgebra of $R / \initial_<(I)$ (see, for example, \cite{SS, KM}).
We provide an example illustrating this functionality below. 

\begin{example}[{\cite[Example~2]{SS}}]
Fix a positive integer $n$. If we let 
\[
S = \QQ[a,b,c,d,u_1, \dots, u_n, v_1, \dots, v_n] / (ad - bc - 1),
\]
and $A \subseteq S$ be the subalgebra generated by $G = \{au_i + bv_i, cu_i + dv_i \}_{1 \le i \le n} \subseteq Q$, then $G \cup \{u_i v_j - u_j v_i \}_{1 \le i < j \le n}$ forms a subalgebra basis for $A$ with respect to \Lex. 
We verify this below when $n = 3$:
\begin{lstlisting}[language = Macaulay2]
i1 : n = 3;
i2 : R = QQ[a,b,c,d,u_1 .. u_n, v_1 .. v_n, MonomialOrder => Lex];
i3 : S = R / ideal(a*d - b*c - 1);
i4 : G = flatten for i from 1 to n list {a*u_i + b*v_i, c*u_i + d*v_i};
i5 : SB = subalgebraBasis G
o5 = | cu_3+dv_3 cu_2+dv_2 cu_1+dv_1 au_3+bv_3 au_2+bv_2 au_1+bv_1 u_2v_3-u_3v_2 u_1v_3-u_3v_1 u_1v_2-u_2v_1 |
             1       9
o5 : Matrix S  <--- S
i6 : isSAGBI SB
o6 = true
\end{lstlisting}
\end{example}




\section{Design and Functionality}\label{sec:design-and-functionality}

\subsection{Data structures and resuming computations}\label{sec:data-structures} We provide two main data structures for working with subalgebras and subalgebra bases: the \subring\ and \SAGBIBasis\ types. 
An instance of the type \subring, which we call a \emph{subring}, represents a subalgebra $A$ of a polynomial ring or a quotient ring.
Subrings are designed to behave similarly to the core data type \Ideal. The other main type \SAGBIBasis\ captures the state of a subalgebra basis computation, making it more similar to a \GroebnerBasis\ computation object. In particular, this type is designed to keep track of the already-computed elements of the subalgebra basis, and we refer to an instance of \SAGBIBasis\ as a \textit{computation object}. 

A subring \texttt{A} is a light-weight object that keeps track of the generators of the algebra $A$, and may be used as an argument for the main functions that compute subalgebra bases;  \texttt{subalgebraBasis} and \texttt{sagbi}. 
A subring also keeps track of the furthest-advanced computation object that has been constructed during a subalgebra basis computation. To access or create a computation object, we use the function \texttt{sagbi}. The generators associated with a computation object are the currently-known elements of a (partial) subalgebra basis. Applying the method function \texttt{isSAGBI} to a subring
checks whether the generators of furthest advanced computation object is a certified subalgebra basis. 

Algebra generators associated with an instance of \subring\ or \SAGBIBasis\ can be recovered using \gens. 
For example, if \texttt{S} is a subring and \texttt{isSAGBI S} returns \true, then a complete subalgebra basis for \texttt{S} has been computed. In this case, the subalgebra basis may be accessed with \gens\ \texttt{sagbi S}, or equivalently \texttt{subalgebraBasis S}.

\begin{example}
Let $A=\QQ [x+y,x^6,y^6] \subseteq \QQ [x,y]$.  The following code computes a partial \GRevLex\ subalgebra basis:
\begin{lstlisting}[language=Macaulay2]
i1 : R = QQ[x,y];
i2 : A = subring {x+y, x^6, y^6}
o2 = subring of R with 3 generators
o2 : Subring
i3 : SB = sagbi(A, Limit => 5)
o3 = Partial SAGBIBasis Computation Object with 1 generators, Limit = 5.
o3 : SAGBIBasis
i4 : isSAGBI SB
o4 = false
\end{lstlisting}
The computation object \texttt{SB} above records the state of a subalgebra basis computation accounting for binomial syzygies of degree up to $5.$
This computation does not compute the S-polynomial $(x+y)^6 - (x^6)$, whose initial term is needed to generate the initial algebra.
Since the subring \texttt{A} has recorded the progress of the last computation, we may resume computation with a higher value for \Limit .
\begin{lstlisting}[language=Macaulay2]
i5 : SB' = sagbi(A, Limit => 100)
o5 = SAGBIBasis Computation Object with 3 generators, Limit = 100.
o5 : SAGBIBasis
i6 : isSAGBI SB'
o6 = true
i7 : gens SB'
o7 = | x+y y6 6x5y+15x4y2+20x3y3+15x2y4+6xy5 |
             1       3
o7 : Matrix R  <--- R

\end{lstlisting}
\end{example}

\begin{remark} 
\subring\ and \SAGBIBasis\ both inherit from \HashTable , and thus both types are immutable.
However, the constructors for either type provide the special key \cache\ used to keep track of the progress of computations.
Any subring has two keys  \generators\ and \texttt{ambientRing}, which store the input generators and the ring containing these generators. 
In contrast, many keys of these objects are reserved for internal use, but may be accessed through accessor functions. Two keys of particular interest are the keys \texttt{SAGBIdata} and \texttt{SAGBIoptions}. The key \texttt{SAGBIoptions} stores the options that are used in subalgebra basis computations. The key \texttt{SAGBIdata} stores user-readable information such as the original generators of the subring, the partially computed subalgebra basis, and whether the subalgebra basis computation is complete.
\end{remark}

\subsection{Computation options}\label{sec:options}

When calling \texttt{sagbi} or \texttt{subalgebraBasis}, several
options for fine-tuning computations may be used. 
These options typically carry over when resuming computations. 
Two exceptions are the options \Limit, which is always taken to be the specified or default value, and \texttt{PrintLevel}, which is described below and controls the verbosity of the computation. However, by using the option \texttt{RenewOptions},
the options may be modified. 
The option \texttt{Recompute} is used to completely restart a subalgebra basis computation. 
We catalogue several other options below, which may be useful in various different settings:





\begin{enumerate}[
    align = left,
    leftmargin=*
    ]
    \item[(\texttt{PrintLevel})] This option takes a non-negative integer and controls the verbosity of the function. For successively higher values, the function will print more data related to the computation.
    \item[(\texttt{SubductionMethod})] This option controls whether subduction is performed by top-level code \quoteoption{"Top"} or by engine-level compiled code \quoteoption{"Engine"}. The engine-level code can be faster for computations that require many subduction steps. One advantage of \quoteoption{"Top"} is that one can view intermediate results (e.g.,~from subduction) if a high enough \texttt{PrintLevel} is used.

    \item[(\texttt{Strategy})] This option controls how the computation object is modified after new generators are added to the subalgebra basis. There are two primary strategies: \quoteoption{"DegreeByDegree"} and \quoteoption{"Incremental"}. The strategy \quoteoption{"DegreeByDegree"} computes a partial Gr\"obner basis for the \textit{reduction ideal}, whereas the strategy \quoteoption{"Incremental"} computes a full Gr\"obner basis. However, the strategy \quoteoption{"DegreeByDegree"} computes the partial Gr\"obner basis from scratch, whereas the strategy \quoteoption{"Incremental"} makes use of previously computed full Gr\"obner basis to speed up the computation of the next full Gr\"obner basis. Both strategies have strengths and weaknesses: \quoteoption{"DegreeByDegree"} is suited to computations where a large number of new generators are added at a particular degree, and \quoteoption{"Incremental"} is well suited to computations where very few generators are added at each degree. The default option is the strategy \quoteoption{"Master"} which heuristically blends between the \quoteoption{"DegreeByDegree"} and \quoteoption{"Incremental"} in order to gain performance benefits from each method.

    \item[(\texttt{AutoSubduceOnPartialCompletion})] This option  controls whether the subalgebra basis generators are subducted against each other. More precisely, if no new subalgebra basis generators are found at a particular degree, then the current subalgebra basis generators are subducted against each other. This produces a reduced set of generators. We intend for this option to be used when the supplied generators are suspected to be a subalgebra basis; reducing the generators speeds up the computation of the subsequent binomial syzygies.
\end{enumerate}

\subsection{Other functionality}\label{sec:other-functionality}

The function \subduction\ takes a collection of polynomials $g_1, \dots, g_s$ and subducts a given polynomial $f\in R$ (or matrix of polynomials) against them. 
When $g_1, \ldots , g_s$ are encoded using a subring \texttt{A} or computation object, then $g_1, \dots, g_s$ are taken to be the generators of \texttt{S} or (partial) subalgebra basis generators of \texttt{A}, respectively.

\begin{example}\label{ex: simple subduction}
Fix the polynomial ring $\bQ[x,y]$ and let $G = \{x^2 + x, y^2 + y \}$. We subduct the polynomial $f = x^2y^2 + x^3y$ against $G$ with respect to \GRevLex\ as follows:
\begin{lstlisting}[language=Macaulay2]
i1 : R = QQ[x,y];
i2 : G = {x^2 + x, y^2 + 1};
i3 : subduction(G, x^2*y^2 + x^3*y)
      3       2
o3 = x y - x*y
o3 : R
\end{lstlisting}
\end{example}

As in Example~\ref{ex:basicSubalgebra},
we may compute normal forms using the operator \texttt{\%}. 
Analogous to its use with ideals and Gr\"obner bases, the syntax \texttt{f \% S} works for
\texttt{S} of class either \subring\ or \SAGBIBasis.
If no complete subalgebra basis for \texttt{S} of class \subring\ is known, then \texttt{f \% S} falls back on an \textit{extrinsic method} using Gr\"{o}bner bases to compute the normal form. On the other hand, if either a complete subalgebra basis has been computed for \texttt{S} or \texttt{S} is of class \SAGBIBasis, then \texttt{f \% S} subducts \texttt{f} against the (partial) subalgebra basis associated with \texttt{S}.

\begin{example}\label{ex:sagbiPercent}
Consider the subalgebra $A = \bC[x+y, xy, xy^2]$ in Example~\ref{ex:degreeLimitedInfiniteSubalgebraBasis}. Let $f = xy^3 + xy^4 + xy^5 + xy^6 \in A$ be a polynomial. We illustrate the different behaviors of \texttt{\%} in the code below.

\begin{lstlisting}[language=Macaulay2]
i1 : R = QQ[x,y];
i2 : A = subring {x+y, x*y, x*y^2};
i3 : SB = sagbi(A, Limit => 5);
i4 : f = x*y^3 + x*y^4 + x*y^5 + x*y^6;
i5 : f % A
o5 = 0
o5 : R
i6 : f % SB
        6      5
o6 = x*y  + x*y
o6 : R
i7 : SB = sagbi(A, Limit => 7);
i8 : f % SB
o8 = 0
o9 : R
\end{lstlisting}

The generators associated with \texttt{SB}, namely $\{x+y, xy, xy^2, xy^3, xy^4, xy^5, xy^6\}$, form a partial subalgebra basis for $A$ . To subduct $f$ against these generators, we may use \texttt{f \% SB}. However, the generators of \texttt{SB} do not form a complete subalgebra basis for \texttt{A}, so in this case the operation \texttt{f \% A} falls back on the extrinsic method.
\end{example}

The function \texttt{groebnerMembershipTest} is an extrinsic membership test for elements of a subring. It can potentially be used in cases where a subring has an intractable subalgebra basis. If a sufficiently large partial subalgebra basis has already been computed, then it is recommended to use the operator \texttt{\%} or the function \texttt{subduction}.


\begin{example}\label{ex: groebnerMembershipExample}
Following from Example~\ref{ex:sagbiPercent}, we test the membership of $f$ in $A$ as follows:
\begin{lstlisting}[language=Macaulay2]
i7 : groebnerMembershipTest(f, A)
o7 = true
\end{lstlisting}
\end{example}

The function \texttt{subringIntersection} computes the intersection of two subrings using an analogous method to that of intersecting ideals via Gr\"obner bases \cite{cox2013ideals}. 
The output of the function is an instance of \texttt{IntersectedSubring}, which is a type that inherits from \subring. 
To check whether the output is guaranteed to be equal to the full intersection of the subrings, we use the function \texttt{isFullIntersection}. 
If the function \texttt{isFullIntersection} returns \texttt{true} then the generators of the intersected subring are a subalgebra basis for the intersection.

\begin{example}
Consider the subrings $A_1 = \QQ[x^2, xy]$ and $A_2 = \QQ[x, y^2]$ of the quotient ring $S = \QQ[x, y]/\langle x^3 + xy^2 + y^3 \rangle$. We compute the intersection $A_1 \cap A_2$ as follows:
\begin{lstlisting}[language=Macaulay2]
i1 : R = QQ[x,y];
i2 : I = ideal(x^3 + x*y^2 + y^3);
i3 : S = R/I;
i4 : A1 = subring {x^2, x*y};
i5 : A2 = subring {x, y^2};
i6 : A = subringIntersection(A1, A2)
o6 = QQ[p_0..p_5], subring of S
o6 : IntersectedSubring
i7 : gens A
o7 = | x2 x2y2+xy3 y4 xy3 y6 xy5 |
             1       6
o7 : Matrix S  <--- S
i8 : isFullIntersection A
o8 = true
\end{lstlisting}
Since \texttt{isFullIntersection A} returns \true , the generators of the intersected subring \texttt{A} are in fact guaranteed to generate the entire ring $A_1 \cap A_2.$ Moreover the generators are a subalgebra basis for the intersection with respect to the default \GRevLex.

\end{example}

If \texttt{isFullIntersection A} returns \false\ for an intersected subring \texttt{A}, then the generators of \texttt{A} may or may not generate the entire intersection. If this happens, then the function \texttt{subringIntersection} may be used with the option \Limit\ set to a value greater than its default of $20$.

\section{Examples}\label{sec: examples}

Subalgebra bases have a number of applications. A short sample of the literature, including several recent papers which use the \SubalgberaBases\ package, follows: \cite{SS, Bossinger2017, CrookDonelan, Gobel, SturmfelsXu, Belotti2022Discrete, ClusterAlgebras, Schubert,Dalbec,breiding2022algebraic,misra2021gaussian}.

\begin{example}\label{ex:screws}
Following \cite[Section~2]{CrookDonelan}, we describe the adjoint action of $\SO(3)$ on its Lie algebra and its ring of invariants. We write $(R, \mathbf{t}) \in \SO(3) \ltimes \RR^3 = \SE(3) $ for an element of the Special Euclidean group of rigid motions of $\RR^3$. To describe the adjoint action of $\SE(3)$ on its Lie algebra $\se(3)$, it is convenient to define, for each $\mathbf{t} = (t_1,t_2,t_3)^T \in \RR^3$, a skew-symmetric matrix 
\[
[\mathbf{t}]_\times  = \begin{pmatrix}
0 & -t_3 & t_2 \cr
t_3 & 0 & -t_1 \cr
-t_2 & t_1 & 0
\end{pmatrix}.
\]
We identify elements of $\se(3)$ with their Pl\"ucker coordinates $(\mathbf{w}, \mathbf{v}) \in \se(3)$. The adjoint action of $\SE(3)$ on $\se(3)$ is given by
\begin{equation}\label{eq:adjoint-action}
(R, \mathbf{t}) \cdot \begin{pmatrix}
\mathbf{w} \cr
\mathbf{v}
\end{pmatrix} =
\begin{pmatrix}
R & 0 \cr
[\mathbf{t}]_\times R & R
\end{pmatrix}
\begin{pmatrix}
\mathbf{w} \cr
\mathbf{v}
\end{pmatrix}
=
\begin{pmatrix}
R \mathbf{w} \cr
[\mathbf{t}]_\times R \mathbf{w} + R \mathbf{v}
\end{pmatrix}.
\end{equation}
Consider the polynomial algebra $A \subseteq \QQ [t_i, w_i, v_i]_{i \in \{1,2,3\}}$ generated by the entries of the matrix in Equation~\eqref{eq:adjoint-action} with the rotation $R$ set to be the $3\times 3$ identity matrix. Explicitly, the subring is given by
\[
A = \QQ[w_1,w_2,w_3, -t_3w_2 + t_2w_3 + v_1, t_3w_1 - t_1w_3 + v_2, -t_2w_1 + t_1w_2 + v_3].
\]
The ring of invariants of $\se(3)$ under the action of the translational subgroup $\RR^3 \triangleleft \SE(3)$ is given by the intersection $A \cap \QQ[w_i, v_i]_{i \in \{1,2,3\}}$. We compute this intersection using a monomial order that eliminates the variables $t_1, t_2, t_3$:
\begin{lstlisting}[language=Macaulay2]
i2 : R = QQ[t_1..t_3, w_1..w_3, v_1..v_3, MonomialOrder=>{Eliminate 3, Lex}];
i3 : SB = sagbi {w_1, w_2, w_3, -t_3*w_2+t_2*w_3+v_1, t_3*w_1-t_1*w_3+v_2, -t_2*w_1+t_1*w_2+v_3};
i4 : isSAGBI SB
o4 = true
i5 : SB' = selectInSubring(1, gens SB)
o5 = | w_3 w_2 w_1 w_1v_1+w_2v_2+w_3v_3 |
             1       4
o5 : Matrix R  <--- R

\end{lstlisting}
It follows from this computation that the algebra of translational invariants is simply $\QQ[w_1, w_2, w_3, \mathbf{w} \cdot \mathbf{v}]$, which confirms the computation in \cite[Section~5.1]{CrookDonelan}. The action above naturally extends to an action on $(\mathbf{w}_1,\mathbf{v}_1, \dots, \mathbf{w}_n, \mathbf{v}_n) \in \se(3)^n$, called a multi-screw. In \texttt{accompanyingCode.m2}, calling \texttt{screwsExample n} will attempt computing a subalgebra basis for the translational invariants for any number $n.$
In particular, it is straightforward to verify the computation when $n=2$, which is used in \cite[Section~5.2]{CrookDonelan} as a step in computing full invariant ring.
\end{example}


Besides invariant-theoretic applications such as Examples~\ref{ex:basicSubalgebra} and~\ref{ex:screws}, computing subalgebra bases is also
a key operation for constructing \emph{toric degenerations}. 
A toric degeneration is a particular type of flat family whose generic fibers are some variety of interest and whose special fiber is a toric variety. 
Various properties of the variety of interest that are preserved under flat limits (e.g.,~dimension and degree) can then be studied via toric degenerations by passing to the toric fiber.
Newton-Okounkov bodies are a closely-related construction, whose applications include counting the number of solutions to classes of polynomial systems 
(see e.g.,~\cite{DuffHeinSottile, BurrSottileWalker, ObatakeWalker,breiding2022algebraic}).
Using subalgebra bases, toric degenerations have been constructed for many families of varieties, including Grassmannians, Flag varieties, Cox-Nagata Rings \cite{SturmfelsXu, Bossinger2017, MohammadiShaw2019MatchingFields, Ollie2022partialflag}.
We conclude with two representative examples.

\begin{example}
Consider the matrix
\[
A = \begin{bmatrix}
1 & 0 & 0 & 1 & -1 & 0 \\
0 & 1 & 0 & -1 & 0 & 1 \\
0 & 0 & 1 & 0 & 1 & -1
\end{bmatrix}
\]
and let $G = \ker(A)$. The columns of $A$ are six points in $\bP^2$, which are the points of intersection of four generic lines.
Following \cite[Example~2.6]{SturmfelsXu}, the Cox-Nagata ring is given by
\[
R^G = \QQ [x_1, \dots, x_6, L_{124}, L_{135}, L_{236}, L_{456}, M_{16}, M_{25}, M_{34}] \subseteq \QQ [x_1, \dots, x_6, y_1, \dots, y_6],
\]
where the polynomials $L_{\bullet}$ and $M_{\bullet}$ are described below. We verify that the generators indeed form a subalgebra basis as follows:
\begin{lstlisting}[language=Macaulay2]
i1  : R = QQ[x_1..x_6, y_1..y_6];
i2  : L124 = y_3*x_5*x_6 + x_3*y_5*x_6 - x_3*x_5*y_6;
i3  : L135 = y_2*x_4*x_6 - x_2*y_4*x_6 + x_2*x_4*y_6;
i4  : L236 = y_1*x_4*x_5 + x_1*y_4*x_5 - x_1*x_4*y_5;
i5  : L456 = y_1*x_2*x_3 + x_1*y_2*x_3 + x_1*x_2*y_3;
i6  : M16 = y_2*x_3*x_4*x_5 + x_2*y_3*x_4*x_5 - x_2*x_3*y_4*x_5 + x_2*x_3*x_4*y_5;
i7  : M25 = y_1*x_3*x_4*x_6 + x_1*y_3*x_4*x_6 + x_1*x_3*y_4*x_6 - x_1*x_3*x_4*y_6;
i8  : M34 = y_1*x_2*x_5*x_6 + x_1*y_2*x_5*x_6 - x_1*x_2*y_5*x_6 + x_1*x_2*x_5*y_6;
i9  : RG = subring {x_1 .. x_6, L124, L135, L236, L456, M16, M25, M34};
i10 : isSAGBI RG
o10 = true
\end{lstlisting}
\end{example}

\begin{example}\label{ex:PluckerPickingUpWhereItLeftOff}
Let $X = (x_{i,j})$ be a $3 \times 6$ matrix of variables and $R = \QQ [x_{i,j}]$ be a ring in those variables. Under the Pl\"ucker embedding, the coordinate ring of the Grassmannian Gr$(3,6)$ is generated by the maximal minors of $X$. However, the maximal minors do not form a subalgebra basis with respect with respect to a weight vector associated to a \textit{hexagonal matching field} \cite{MohammadiShaw2019MatchingFields,Ollie2022partialflag}. 
\begin{lstlisting}[language=Macaulay2]
i1 : R = QQ[x_(1,1)..x_(3,6), MonomialOrder => {Weights => {0,0,0,0,0,0,
0,15,3,12,9,6,0,7,14,21,28,35}}];
i2 : X = transpose genericMatrix(R, 6, 3);
             3       6
o2 : Matrix R  <--- R
i3 : A = subring for s in subsets(6, 3) list det X_s;
i4 : SB = sagbi(A, Limit => 100, SubductionMethod => "Engine")
o4 = Partial SAGBIBasis Computation Object with 21 generators, Limit = 100.
o4 : SAGBIBasis
i5 : isSAGBI SB
o5 = true
\end{lstlisting}
Note that the first $20$ generators of \texttt{SB} are the maximal minors of $X$, or, in other words, the original generators of \texttt{A}. However, a twenty-first generator is required for a complete subalgebra basis.
\end{example}

\section*{Acknowledgements}
We wish to thank the organizers of two virtual \Macaulay\ workshops held in 2020 virtually ``at'' Cleveland State University and the University of Warwick. 
These workshops provided first introductions for several of the authors. 
We also thank several working group participants who we also wish to thank for their help during early stages of re-developing this package.

This article has been co-authored by an employee of National Technology \& Engineering Solutions of Sandia, LLC under Contract No. DE-NA0003525 with the U.S. Department of Energy (DOE). The employee co-owns right, title and interest in and to the article and is responsible for its contents. The United States Government retains and the publisher, by accepting the article for publication, acknowledges that the United States Government retains a non-exclusive, paid-up, irrevocable, world-wide license to publish or reproduce the published form of this article or allow others to do so, for United States Government purposes. The DOE will provide public access to these results of federally sponsored research in accordance with the DOE Public Access Plan \url{https://www.energy.gov/downloads/doe-public-access-plan}.
\bibliographystyle{plain}
\bibliography{references}

\end{document}